\documentclass{amsart}

\usepackage{amsthm,amssymb,latexsym,amsmath}
\usepackage{mathrsfs}
\usepackage{graphicx}
\input xy
\xyoption{all}
\usepackage[all]{xy}

\newcommand{\F}{ {\mathcal F}}
\renewcommand{\sec}{\mathbb{S}ec}

\newtheorem{nada}{Nada}[section]
\newtheorem{definition}[nada]{Definition}

\newtheorem{theorem}[nada]{Theorem}
\newtheorem{lemma}[nada]{Lemma}

\newcommand{\bc}{\begin{center}}
\newcommand{\ec}{\end{center}}

\begin{document}

\title{ A  singular Darboux  type theorem and  non-integrable  projective distributions of degree one  }
\hyphenation{ho-mo-lo-gi-cal} \hyphenation{fo-lia-tion}

\dedicatory{\it To Omegar Calvo-Andrade on the occasion of his 60th birthday}

\begin{abstract}  
We prove a singular Darboux  type theorem for   homogeneous   polynomial  closed  $2$-forms of  degree one on $\mathbb{C}^n$.
As application, we classify non-integrable codimension one  distributions, of degree one, and  arbitrary classes  on  projective spaces. 
\end{abstract}

\author{Maur\'icio Corr\^ea}
\address{Maur\'icio Corr\^ea \\ ICEx - UFMG \\
Departamento de Matem\'atica \\
Av. Ant\^onio Carlos 6627 \\
30123-970 Belo Horizonte MG, Brazil } \email{mauriciojr@ufmg.br}
\author{Vin\'icius Soares dos Reis}
\address{ Vin\'icius Soares dos Reis\\ ICEB - UFOP\\ Departamento de Matem\'atica\\ Campus Universitário Morro do Cruzeiro, 35400-000,  Ouro Preto MG, Brazil,}
\email{dosreis.vinicius@gmail.com}
\thanks{ }

\maketitle

\section{Introduction}
The classical Darboux  Theorem  states that  if  $\omega$ is  a   closed  non-singular  holomorphic   $2$-form on  $\mathbb{C}^n$
which satisfies $\omega^{k} \neq 0$  and  $\omega^{k+1}=0$, then there exists a coordinate system  $(x_{1},\dots, x_{k}, y_{1}, \dots, y_{k},z)\in \mathbb{C}^{n}$, with $z=(z_{2k+1},\dots,z_{n})$ such that
$$\omega =  dx_{1}\wedge dy_{1} +\cdots +dx_{k}\wedge dy_{k}.$$

In this work we prove  the following  version of  Darboux   Theorem for degree one    homogeneous   polynomial differential $2$-forms     on $\mathbb{C}^n$.

\begin{theorem}\label{TeoDarb}
Let $ \omega $ be a homogeneous polynomial closed  $2$-form of degree one on  $\mathbb{C}^{n}$ such that $\omega^{k}   \not\equiv 0$  and  $\omega^{k+1} \equiv 0$. Then, there exists a coordinate system $(x_{1},\dots, x_{k}, y_{1}, \dots, y_{k},z) \in \mathbb{C}^{n}$, with $z=(z_{2k+1},\dots,z_{n})$ 
such that $\omega$ reduces to one of the following normal forms:
\begin{enumerate}
\item 
$
\omega=\pi^*\eta, 
$
where  $\eta$  is a   closed  homogeneous $2$-form  of degree one  on  $\mathbb{C}^{2k}$  and   $\pi(x_{1},\dots, x_{k}, y_{1}, \dots, y_{k},z)=(x_{1},\dots, x_{k}, y_{1}, \dots, y_{k}) $ denotes the canonical linear projection. 
\item
$\omega =   \pi^* \vartheta +dt_{1}\wedge dh_{1} +\cdots +dt_{k}\wedge dh_{k},$
where  $\vartheta $  is a   closed  homogeneous $2$-form  of degree one  on  $\mathbb{C}^{2k}$, $t_{1},\dots, t_{k}$ are linear  polynomials  in the variables $(x_{1},\dots, x_{k}, y_{1}, \dots, y_{k})$ and $h_{1}, \dots, h_{k}$ are quadratic  polynomials.
\end{enumerate}
\end{theorem}

This  result is   a version,  for non locally decomposable $2$-forms,   of a  Theorem due to A. Medeiros \cite[Theorem A]{De Medeiros}. More precisely,  
 Medeiros  proved the following:  if   $ \omega $  is a locally decomposable  closed 
homogeneous $q$-form of degree  one   on  $\mathbb{C}^{n}$, which is not a linear pull-back, then 
$$
\omega=dt_{1}\wedge \dots \wedge dt_{q-1}\wedge dh,
$$
where $t_{1},\dots, t_{q-1}$ are linear  polynomials and $h$ is a quadratic polynomial. 
 Let us show how we can recover the  Medeiros's result  in the case  $q=2$  and $k=1$.   By Theorem \ref{TeoDarb}, if $\omega$ is not a linear pull-back   we can write  
$$\omega =   g(x , y)dx\wedge dy+(adx+bdy)\wedge dh, $$
 where $g(x , y)$ is linear and  $h$ is  quadratic. 
 If $b=0$, then  $\omega =    ( -g(x , y)dy-adh)\wedge dx$.   Otherwise, 
we can write
 $\omega =   (b^{-1}g(x , y)dx -dh)\wedge (adx+bdy) $. That is,  we can always  write 
   $$   \omega =   \alpha\wedge dt$$
   for some quadratic $2$-form $\alpha$ and some linear polynomial $t$. 
Now,   $ 0=d\omega =  d \alpha\wedge dt$ implies that $d\alpha=d\ell\wedge dt=d(\ell dt)$, for some linear polynomial $\ell$.  Therefore, there exist a  polynomial $q$ of degree $2$   such that
$\alpha- \ell dt=dq$. Therefore, 
  $$   \omega =   \alpha\wedge dt=(dq+\ell  dt)\wedge dt=dq\wedge dt. $$

In \cite{Jouanolou}  J-P   Jouanolou  classified  codimension one holomorphic  foliations of degree one on $\mathbb{P}^{n}$. More precisely, he  showed  that if $\F$ is such foliation, then we
are in one of following cases: 
\begin{enumerate}
\item[i)] $\mathcal{F}$ is defined a dominant rational map $\mathbb{P}^{n}\dashrightarrow \mathbb{P}(1,2)$ with irreducible general fiber determined  by  a linear  polynomial  and one quadratic polynomial; or
\item[ii)] $\mathcal{F}$ is the linear pull back of a foliation of induced by a global holomorphic vector field on $\mathbb{P}^{2}.$
\end{enumerate}

Loray, Pereira and Touzet in [\ref{LPT} , Theorem 6] generalized the  Theorem of  Jouanolou  to     holomorphic foliation   with codimension $q\geq 1$   and  degree one by using  the 
  above cited result due to Medeiros  \cite[Theorem A]{De Medeiros}. They showed that  
 if $\mathcal{F}$  is a foliation of degree  one  and codimension $q$ on $\mathbb{P}^{n}$, then  we
are in one of following cases 
\begin{enumerate}
\item[i)] $\mathcal{F}$ is defined a dominant rational map $\mathbb{P}^{n}\dashrightarrow \mathbb{P}(1^{q},2)$ with irreducible general fiber determined  by $q$  linear  polynomials and one quadratic form; or
\item[ii)] $\mathcal{F}$ is the linear pull back of a foliation of induced by a global holomorphic vector field on $\mathbb{P}^{q+1}.$
\end{enumerate}

 Let $\F $ be a codimension one distribution on a complex  manifold  $X$, and consider the  associated 
line bundle $L_\F:=\det(T_X/T\F)$, where $T\F$ denotes its tangent sheaf.  
The distribution $\F$ corresponds to a unique (up to scaling) twisted $1$-form 
$$\omega_{\F}\in H^0(X,\Omega^1_X\otimes L_\F)$$ 
 non vanishing in codimension one. 
For every integer $i\geq 0$, there is a well defined twisted $(2i+1)$-form 
$$
\omega_{\F}\wedge (d\omega_{\F})^i\ \in \ H^0\Big(X,\Omega^{2i+1}_X\otimes L_\F^{\otimes (i+1)}\Big).
$$
The \emph{class} of  $\F$ is the unique non negative  integer $k=k({\F})$ such that 
$$
\omega\wedge (d\omega)^k  \not\equiv 0\  \ \text{ and } \ \ \omega\wedge (d\omega)^{k+1} \equiv0.
$$
By Frobenius theorem, a codimension one distribution is a foliation if and only if $k({\F})=0$.

We shall use our Darboux type theorem  in order to classify  non-integrable distributions  on  $\mathbb{P}^{n}$ of degree one  and arbitrary class.

\begin{theorem}\label{classif}
 Let $\mathcal{F}$ be a distribution on  $\mathbb{P}^{n}$  of degree  one  and class $k\geq 1$. Then  we
are in one of following cases: 
\begin{enumerate}
\item[i)] there is a rational linear map $\rho: \mathbb{P}^n  \dashrightarrow \mathbb{P}^{2k+1}$ and a distribution $\mathcal{G}$ of degree  one  on  $ \mathbb{P}^{2k+1}$ such that $\mathcal{F}= \rho^*\mathcal{G}$.
\item[ii)]  there is a rational map   $\xi: \mathbb{P}^n \dashrightarrow  \mathbb{P}(1^{k+1},2^{k+1})$ 
determined by $k+1$ linear polynomials and  $k+1$ quadratic  polynomials, and  a  rational linear map $\rho: \mathbb{P}^n  \dashrightarrow \mathbb{P}^{2k+1}$  such that  $\mathcal{F}$ is induced by $ \rho^*\alpha +\xi^*\theta_0$, where $\alpha \in  H^0(  \mathbb{P}^{2k+1}, \Omega_{ \mathbb{P}^{2k+1}}^1(3))$  and   $\theta_{0}=   \sum_{i}(u_{i}dw_{i} - 2w_{i}du_{i}) $  is the canonical contact distribution on  $\mathbb{P}(1^{k+1},2^{k+1})$.
\end{enumerate}
\end{theorem}

We observe that  $\alpha \in  H^0(  \mathbb{P}^{2k+1}, \Omega_{ \mathbb{P}^{2k+1}}^1(3))$ can be zero and in this case we have that  
 $\mathcal{F}= \xi^*\mathcal{G}_0$, where $ \mathcal{G}_0$ is the canonical contact distribution on  $\mathbb{P}(1^{k+1},2^{k+1})$. 
In  \cite{CCF}   the authors have showed that under generic conditions    a  non-integrable distribution on $\mathbb{P}^n$  is this type. 

Now, consider  $\mathcal{D}(d;k;n) \subset \mathbb{P} H^0(X,\Omega^1_{\mathbb{P}^n}(d+2))$ the space of distributions of codimension one  on  $\mathbb{P}^{n}$, of degree $d$ and class $k$.  We can see that the spaces $\mathcal{D}(d;k;n) $ are algebraic subvarieties of $ \mathbb{P} H^0(X,\Omega^1_{\mathbb{P}^n}(d+2))$. 
It is well known that the space  $\mathcal{D}(0;0;n)$ of degree zero foliation is  the Grassmannian of lines in $\mathbb{P}^n$.  For more details about the spaces of foliations  $\mathcal{D}(d;0;n)$  see  \cite{Omegar} and \cite{Lins Neto}   references therein.  

In \cite{Carolina} Ara\'ujo, Corr\^ea and Massarenti studied  in particular the geometry of spaces $\mathcal{D}(0;k;n)$.  More precisely, the authors in \cite{Carolina}  showed the following:

Let $D_k\subseteq\mathbb{P}(H^0(\mathbb{P}^n,\Omega^1_{\mathbb{P}^n}(2)))$ be the variety parametrizing codimension one distributions on 
$\mathbb{P}^n = \mathbb{P}(\mathbb{C}^{n+1})$ of class $\leq k$ and degree zero. 
Identify $H^0(\mathbb{P}^n,\Omega^1_{\mathbb{P}^n}(2))$ with  $\bigwedge^2\mathbb{C}^{n+1}$.
Then $D_k = \sec_{k+1}(\mathbb{G}(1,n))$ and the stratification 
$$D_0\subseteq D_1\subseteq ...\subseteq D_{k-1}\subseteq ...\subseteq \mathbb{P}(H^0(\mathbb{P}^n,\Omega^1_{\mathbb{P}^n}(2)))$$
corresponds to the natural stratification
$$\mathbb{G}(1,n)\subseteq \sec_2(\mathbb{G}(1,n))\subseteq ...\subseteq \sec_k(\mathbb{G}(1,n))\subseteq ...\subseteq \mathbb{P}(\bigwedge^2\mathbb{C}^{n+1}),$$
where   $\sec_j(\mathbb{G}(1,n))$ is the 
  \emph{$i$-secant variety} of       the  Grassmannian $\mathbb{G}(1,n)$  of lines in $\mathbb{P}^n$.
  
In \cite{CCJ16} the  Calvo-Andrade, Corr\^ea and Jardim  have studied   codimension one distributions of class one and low degrees on $\mathbb{P}^3$ describing their moduli spaces  
in terms of moduli spaces of stables sheaves. 
In order to describe the geometry of the spaces $\mathcal{D}(1;k;n)$ we   believe   that  Theorem \ref{classif}  might be useful. 
We  will not consider this problem  in this  work.

  \subsection*{Acknowledgments}
We are  grateful to Marcio Soares,  Bruno Sc\'ardua, Arturo Fernandez-Perez and  A. M. Rodr\'iguez for interesting conversations.
The first named author was partially supported  by CNPq, CAPES and  FAPEMIG.   The first named author is grateful to  University  of Oxford for hospitality. 
Finally, we would like to thank the  anonymous referee  for valuable remarks and suggestions.

\section{Polynomial differential $r$-forms}

Consider the exterior algebra of polynomials differential  $r$-forms in $\mathbb{C}^{n}$ given by $$\Omega^{r}(n) := \wedge^r(\mathbb{C}^{n}) \otimes \mathbb{C}[z],$$
where $\mathbb{C}[z]:=\mathbb{C}[z_1,\dots,z_n]$.
Let $S_{d}$ be the subspace of $\mathbb{C}[z]$ of polynomials of degree $\leq d$. The algebra  $\Omega^{r}(n)$ is naturally graduated:
$$\Omega^{r}(n) = \displaystyle\bigoplus_{d \in \mathbb{N}}\Omega_{d}^{r}(n),$$
where $\Omega_{d}^{r}(n)=  \wedge^r(\mathbb{C}^{n}) \otimes S_{d}.$ 

We will  denote  the module generated by the differentials $dz_{i_1},\dots, dz_{i_r}$, with $i_1< \dots < j_r$,  by $\langle  dz_{i_1},\dots, dz_{i_r} \rangle$. That is, if  
$\Omega\in\langle dz_{i_1},\dots, dz_{i_r}\rangle$  is a  polynomial  differential $1$-form, then $\Omega=\sum_j f_{i_j}dz_{i_j} $
 with $f_{i_j}\in \mathbb{C}[ z] $, for all $j=1,\dots,r.$

Now, consider  a polynomial $r$-form $\omega$:
$$\omega = \displaystyle\sum_{1 \leq i_{1} < ... < i_{r} \leq n}P_{i_{1} ... i_{r}}dz_{i_{1}}\wedge...\wedge dz_{i_{r}}.$$
The degree of $\omega$ is defined by deg$(\omega)=$ max\{deg$(P_{i_{1},...,i_{r}}), 1\leq i_{1} < ... < i_{r} \leq n\},$ and, if $\omega \in \Omega_{d}^{r}(n)$, then $P_{i_{1},...,i_{r}} \in S_{d}.$

Consider  the radial vector field on $\mathbb{C}^{n}$ which is given by 
$$
R=x_1 \frac{\partial}{\partial x_1} + \cdots + x_n
\frac{\partial}{\partial x_n}.
$$

\indent In order to proof the theorem \ref{classif}  we  shall use the following Jouanoulou's Lemma.
\begin{lemma}{\textbf{(Jouanoulou's Lemma) [\ref{jo}, Lemme 1.2, pp. 3]}}\label{170101}
If $\eta$ is a homogeneous polynomial differential  $q$-form of degree $s$, then
$$i_{R}d\eta+d(i_{R}\eta) = (q+s)\eta$$
where $R$ is the radial vector field and $i_{R}$ denotes the interior product or contraction with $R$.
\end{lemma}
See \cite{Correa Maza Soares} and \cite{CM} for more details about   polynomial differential systems.

\section{Codimension one holomorphic distributions} \label{section:distributions}

Let $X$ be a complex manifold of dimension $n$.

\begin{definition}
A  holomorphic distribution $\F$ on $X$ of codimension one is a  nonzero subsheaf $T\F\subset T_X$, of rank $(n-1)$, which is saturated, i.e. such that the quotient $T_X/T\F$ is torsion-free. The sheaf $T\F$ is called by the \emph{tangent sheaf} of $\F$. 
The \emph{normal sheaf} of $\F$ is the   sheaf $N_\F:=(T_X/T\F) $. We denote its determinant by $L_{\F}=\det(N_\F)$.
The singular locus of $\F$ is the locus $Sing(\F)$ where $N_\F$ fails to be locally free.

\end{definition}

Let $\F$ be a codimension  one  distribution on $X$.
The $(n-1)$-th wedge product of the inclusion $N^*_\F\subset \Omega^1_X$ gives rise to a  twisted $1$-form  
$$\omega_{\F}\in  H^0(X,\Omega^1_X\otimes  L_\F)$$ non vanishing in codimension $1$. 
We have that  $T\F$ is the  kernel
of the morphism   given by the contraction with $\omega_{\F}$. 

For every integer $i\geq 0$, there is a well defined twisted $(2i+1)$-form 
$$
\omega_{\F}\wedge (d\omega_{\F})^i\ \in \ H^0\big(X,\Omega^{2i+1}_X\otimes L_\F^{\otimes (i+1)}\big).
$$
The \emph{class} of  $\F$ is the unique non negative  integer $k=k({\F})\in\big\{0, \cdots, \lfloor\frac{n-1}{2}\rfloor\big\}$ such that 
$$
\omega_{\F}\wedge (d\omega_{\F})^k\neq0\  \ \text{ and } \ \ \omega_{\F}\wedge (d\omega_{\F})^{k+1}=0.
$$
See \cite{Godbillon} and  \cite{Bryant}.
 
\subsection{Holomorphic projective distributions}
A  holomorphic   distribution   on  a complex
projective space will be called by   holomorphic projective distribution.

Let  $H^0(\mathbb{P}^n,\Omega^{1}_{\mathbb{P}^n}\otimes L_{\F})$ be a   twisted $1$-form induced by a 
codimension  one  distribution  $\F$ on   a complex
projective space $\mathbb{P}^n.$

If $i: \mathbb P^1 \to \mathbb P^n$ is a generic linear immersion, 
then $i^* \omega_{\mathcal{I}} \in H^0(\mathbb P^1,
\Omega^1_{\mathbb P^1} \otimes  L_{\F})$ is a section of a line
bundle, and its divisor of zeros reflects the tangencies between
$\F$ and $i(\mathbb P^1)$. The \emph{degree} of
$\F$ is, by definition, the degree of such a tangency
divisor. Set $d:=\deg(\F)$. Since $\Omega^1_{\mathbb
P^1}\otimes \mathcal L = \mathcal O_{\mathbb P^1}( \deg(  L_{\F})
-2)$, one concludes that $  L_{\F} = \mathcal O_{\mathbb
P^n}(d+ 2)$. That is, a 
codimension  one  distribution  $\F$ on   a complex
projective space $\mathbb{P}^n$  of degree $d$ induces  a global section   $\omega\in
H^0(\mathbb{P}^n,\Omega^{1}_{\mathbb{P}^n}(d+2)).$

Besides, the Euler sequence implies that a section $\omega$ of
$\Omega^1_{\mathbb P^n} ( d + 2 )$ can be thought of as a
polynomial $k$-form on $\mathbb{C}^{n+1}$ with homogeneous
coefficients of degree $d + 1$, which we will still denote by
$\omega$, satisfying
\begin{equation}
\label{equirw} i_R \omega = 0
\end{equation}
where
$$
R=x_0 \frac{\partial}{\partial x_0} + \cdots + x_n
\frac{\partial}{\partial x_n}
$$
is the radial vector field.

\section{Proof of Theorem \ref{TeoDarb}}

Let $\omega$ be a homogeneous  polynomial differential $2$-form of degree $1$ such that,
\begin{equation*}\label{090101}
\omega^{k}   \not\equiv 0 ;\hspace{1 cm}  \omega^{k+1} \equiv  0.
\end{equation*}
Taking $p_{0}\in\mathbb{C}^{n}\setminus Sing (\omega)$, it follows from  the classical  Darboux theorem that  we can find  a neighborhood of $x_{0}$ and   $1$-forms  $\theta_{i}, \alpha_{i}$ $i\in\{1,\cdots k\}$ such that
\begin{equation}
\omega = \theta_{1}\wedge\alpha_{1}+\cdots +\theta_{k}\wedge\alpha_{k}.
\end{equation}
We can see $\omega$ as a holomorphic map $\omega: \mathbb{C}^n\to \wedge^2(\mathbb{C}^n)$ and since   $\omega$ is linear 
$\omega=\omega'(p_0)$,  where $\omega'(p_0)$ is the derivative
of $\omega$ at $p_0$. Therefore, we have that 
$$
\omega=\omega'(p_0)= \theta'_{1}(p_{0})\wedge\alpha_{1}(p_{0})-\theta_{1}(p_{0})\wedge\alpha_{1}'(p_{0})+\cdots+ \theta'_{k}(p_{0})\wedge\alpha_{k}(p_{0})-\theta_{k}(p_{0})\wedge\alpha_{k}'(p_{0}).
$$
Now, we can define
\begin{eqnarray*}
 \alpha_{1}(p_{0}) = dx_{1}, &\dots,& \alpha_{k}(p_{0})=dx_{k}, \\
 \theta_{1}(p_{0}) = dy_{1}, &\dots,& \theta_{k}(p_{0})=dy_{k}, \\
 \alpha'_{1}(p_{0}) = \pi_{1}, &\dots,& \alpha'_{k}(p_{0})=\pi_{k}, \\
 \theta'_{1}(p_{0}) = \eta_{1}, &\dots,& \theta'_{k}(p_{0})=\eta_{k}.
 \end{eqnarray*}
Thus, we get that 
\begin{equation}\label{100101}
\omega = \eta_{1}\wedge dx_{1}+\cdots+ \eta_{k}\wedge dx_{k}+\cdots +\pi_{1}\wedge dy_{1}+\cdots+ \pi_{k}\wedge dy_{k}.
\end{equation}
We write 
\begin{eqnarray}\label{090105}
  \eta_{1} &=& l_{12}dx_{2}+\cdots +l_{1k}dx_{k} +m_{11}dy_{1}+\cdots+ m_{1k}dy_{k}+\overline{\eta_{1}},\nonumber \\
  \eta_{2} &=& l_{21}dx_{1}+\cdots +l_{2k}dx_{k} +m_{21}dy_{1}+\cdots+ m_{2k}dy_{k}+\overline{\eta_{2}}, \nonumber \\
  \vdots &=& \vdots \nonumber \\
  \eta_{k} &=& l_{k1}dx_{1}+\cdots +l_{k-1k}dx_{k-1} +m_{k1}dy_{1}+\cdots+ m_{kk}dy_{k}+\overline{\eta_{k}},  \\
 \pi_{1} &=& g_{12}dy_{2}+\cdots +g_{1k}dy_{k} +h_{11}dx_{1}+\cdots+ h_{1k}dx_{k}+\overline{\pi_{1}},\nonumber \\
  \pi_{2} &=& g_{21}dy_{1}+\cdots +g_{2k}dy_{k} +h_{21}dx_{1}+\cdots+ h_{2k}dx_{k}+\overline{\pi_{2}},\nonumber  \\
  \vdots &=& \vdots\nonumber  \\
  \pi_{k} &=& g_{k1}dy_{1}+\cdots +g_{k-1k}dy_{k-1} +h_{k1}dx_{1}+\cdots+ h_{kk}dx_{k}+\overline{\pi_{k}}\nonumber ,
\end{eqnarray}
where $\overline{\eta_{i}},\overline{\pi_{i}} \in \langle dz_{2k+1},\dots,dz_{n}\rangle$ and $g_{ij}, h_{ij}, l_{ij},m_{ij}$ are linear functions for $i, j \in \{1,\dots, k \}$. Therefore,  we have the following expression for $\omega$:

$$
\omega =\sum_{i<j}(l_{ji} - l_{ij})dx_{i}\wedge dx_{j} + \sum_{i<j}(g_{ji}-g_{ij})dy_{i}\wedge dy_{j} + \sum_{i,j}(h_{ji}-m_{ij})dx_{i}\wedge dy_{j} +
$$
$$
 +\sum_{i}\overline{\eta_{i}}\wedge dx_{i} + \sum_{i}\overline{\pi_{i}}\wedge dy_{i}.
$$
By  using the hypothesis $\omega^{k+1} = 0$, we have the following identity:

$$0 = \left[\sum_{i,j} u_{ij}\overline{\pi_{i}}\wedge \overline{\pi_{j}} + \sum_{i,j} v_{ij}\overline{\eta_{i}}\wedge \overline{\eta_{j}} +
+ \sum_{i,j} w_{ij}\overline{\eta_{i}}\wedge \overline{\pi_{j}}\right]\wedge dx_{1}\wedge \cdots  \wedge dx_{k}\wedge  dy_{1}\wedge   \cdots \wedge dy_{k} +
$$
$$
+  \sum_{i,j} \overline{\eta_{i}}\wedge\overline{\eta_{j}}\wedge\overline{\pi_{2}}\wedge\cdots  \wedge\overline{\pi_{k}}\wedge dx_{i}\wedge dx_{j}\wedge \cdots  \wedge dy_{1}\wedge\cdots \wedge dy_{k}+\cdots +
$$
$$
+\sum_{i} \overline{\eta_{i}}\wedge\overline{\pi_{1}}\wedge\cdots \wedge\overline{\pi_{k}}\wedge dx_{i}\wedge dy_{1}\wedge \cdots \wedge dy_{k}+\sum_{i} \overline{\eta_{1}}\wedge \cdots \wedge\overline{\eta_{k}}\wedge\overline{\pi_{i}}\wedge dx_{1}\wedge \cdots \wedge dx_{k}\wedge dy_{i},$$
where $u_{ij}$ and $v_{ij}$ are products of $(l_{ji}-l_{ij})$ with $(g_{ji}-g_{ij})$ and they have degree $k-2$. 

By taking the  wedge  product   with    $dx_{1}\wedge\cdots \wedge \widehat{dx_{i}}\wedge \cdots \wedge dx_{k}$ we obtain
\begin{equation*}
\overline{\eta_{i}}\wedge\overline{\pi_{1}}\wedge\cdots \wedge\overline{\pi_{k}}\wedge dx_{1}\wedge\cdots \wedge dx_{i}\wedge \cdots \wedge dx_{k}\wedge dy_{1}\wedge \cdots \wedge dy_{k}=0,\hspace{0.4 cm} i\in\{1,\dots, k\}.
\end{equation*}
Analogously, if we take wedge  product   with $dy_{1}\wedge\cdots  \wedge \widehat{dy_{i}}\wedge \cdots \wedge dy_{k}$  we also have 
\begin{equation*}
\overline{\pi_{i}}\wedge\overline{\eta_{1}}\wedge \cdots \wedge \overline{\eta_{k}} \wedge dx_{1}\wedge\cdots \wedge dx_{i}\wedge \cdots \wedge dx_{k}\wedge dy_{1}\wedge \cdots \wedge dy_{k}=0,\hspace{0.4 cm} i\in\{1,\dots, k\}.
\end{equation*}

Since  $\overline{\eta_{i}}, \overline{\pi_{i}} \in \langle dz_{2k+1},\dots,dz_{n}\rangle$  we obtain 
\begin{equation*}
\overline{\eta_{i}}\wedge\overline{\pi_{1}}\wedge \cdots \wedge \overline{\pi_{k}} = 0, \hspace{0.4 cm} i\in\{1,\dots, k\},
\end{equation*}
and 
\begin{equation*}
\overline{\pi_{i}}\wedge\overline{\eta_{1}}\wedge \cdots \wedge \overline{\eta_{k}} = 0, \hspace{0.4 cm} i\in\{1,\dots, k\}.
\end{equation*}
Then, we get that 
\begin{equation*}
\overline{\eta_{i}} = a^{i}_{1}\overline{\pi_{1}}+\cdots+a^{i}_{k}\overline{\pi_{k}},\hspace{0.2 cm}a^{i}_{j}\in \mathbb{C},   \ \forall \  i,j\in\{1,\dots,k\} .
\end{equation*}
Therefore, if $\overline{\pi_{i}}= 0$ for all $i,j\in \{1,\dots,k\}$, we  have that 
\begin{equation}
\omega =\sum_{i<j}f_{ij}dx_{i}\wedge dx_{j} + \sum_{i<j}r_{ij}dy_{i}\wedge dy_{j} + \sum_{i,j}s_{ij}dx_{i}\wedge dy_{j},
\end{equation}
where $f_{ij}, r_{ij}$ and $s_{ij}$ are linear functions. Moreover, since 
  $d\omega = 0$, we conclude that in fact 
$$f_{ij}, r_{ij}, s_{ij}\in\mathbb{C}[x_{1},\dots,x_{k},y_{1},\dots,y_{k}].$$
This proves $(1)$.

\indent Suppose   that $\overline{\pi_{i}}\neq 0$ for some $i\in\{1,\dots,k\}$.  Since  $d\omega=0$, we conclude  from (\ref{100101})  that 
\begin{equation}\label{090103}
0=d\eta_{1}\wedge dx_{1}+\cdots +d\eta_{k}\wedge dx_{k}+d\pi_{1}\wedge dy_{1} + \cdots + d\pi_{k}\wedge dy_{k}.
\end{equation}
By taking  wedge  product  of  (\ref{090103})  with  $dx_{2}\wedge \cdots \wedge dx_{k}\wedge dy_{1}\wedge \cdots \wedge dy_{k}$ we get that 
\begin{equation*}
d\eta_{1}\wedge dx_{1}\wedge dx_{2}\wedge \cdots  \wedge dx_{k}\wedge dy_{1}\wedge \cdots \wedge dy_{k}=0.
\end{equation*}
Then 
\begin{equation*}
d\eta_{1}=\sigma_{1}\wedge dx_{1}+ \cdots \sigma_{k}\wedge dx_{k}+\xi_{1}\wedge dy_{1}+\cdots+ \xi_{k}\wedge dy_{k},
\end{equation*}
where $\sigma_{i}$ and $\xi_{i}$ are constants $1$-forms, for all $i$. Then, there is  a linear $1$-form $$\beta_{1}\in\langle dx_{1},\dots, dx_{k}, dy_{1},\dots, dy_{k}\rangle,$$ such that $d\eta_{1} = d\beta_{1}$.  Since  $d(\eta_{1} - \beta_{1}) =0$ there exists a quadratic function $f^{1}_{1}$ in $\mathbb{C}^{n}$ such that
\begin{equation*}
\eta_{1}= \beta_{1} + df^{1}_{1}.
\end{equation*}
\indent  By an analogous argument, we can find   quadratic polynomials  $f_{i}^{1}, f_{i}^{2}$   and linear $1$-forms forms $\beta_{i}, \mu_{i}\in\langle dx_{1},\dots, dx_{k}, dy_{1},\dots, dy_{k}\rangle$ such that
$$
\eta_{i}= \beta_{i} + df^{i}_{1},\ \  \pi_{i}= \mu_{i} + df^{i}_{1}, \ \forall i=1,\dots,k.
$$
Thus, we obtain that 
\begin{equation}\label{090106}
\omega = (\beta_{1}+df_{1}^{1})\wedge dx_{1} + \cdots +(\beta_{k}+df_{k}^{1})\wedge dx_{k}+(\mu_{1}+df_{1}^{2})\wedge dy_{1} +\cdots +(\mu_{k}+df_{k}^{2})\wedge dy_{k}.
\end{equation}
Note that,
$$\pi_{i}-\overline{\pi_{i}} = g_{i1}dy_{1} +\cdots + g_{(i-1)k}dy_{k} + h_{i1}dx_{1} + \cdots +h_{ik}dx_{k},\hspace{0.2 cm} \forall i\in\{1,\dots,k\}.$$
Now, we define 
$$\delta_{i} = g_{i1}dy_{1} +\cdots + g_{(i-1)k}dy_{k} + h_{i1}dx_{1} + \cdots +h_{ik}dx_{k}-\mu_{i},\hspace{0.2 cm} \forall i\in\{1,\dots,k\}.$$
Observe that  $\delta _{i} \in \langle dx_{1}, \dots, dx_{k},dy_{1},\dots, dy_{k}\rangle, \hspace{0.2 cm} \forall i\in\{1,\dots,k\}.$
Since  $$\pi_{i}=\mu_{i}+df_{i}^{2},$$  we have 
\begin{equation}\label{090104}
\overline{\pi_{i}} = df^{2}_{i} - \delta_{i}.
\end{equation}

Substituting (\ref{090104}) in  (\ref{090105}) we obtain 
\begin{center}
$\eta_{i} = l_{i1}dx_{1}+\cdots +l_{i(i-1)}dx_{(i-1)}+l_{i(i+1)}dx_{(i+1)} +\cdots +l_{ik}dx_{k} +m_{i1}dy_{1}+\cdots+ + m_{ik}dy_{k}+a^{i}_{1}(df^{2}_{1} - \delta_{1})+\cdots+a^{i}_{k}(df^{2}_{k} - \delta_{k})$.
\end{center}
Define 
\begin{center}
$\gamma_{i} = l_{i1}dx_{1}+\cdots +l_{i(i-1)}dx_{(i-1)}+l_{i(i+1)}dx_{(i+1)} +\cdots +l_{ik}dx_{k} +m_{i1}dy_{1}+\cdots+ + m_{ik}dy_{k} -a_{1}^{i}\delta_{1}-\dots-a_{k}^{i}\delta_{k} \in \langle dx_{1}, \dots, dx_{k},dy_{1},\dots, dy_{k}\rangle$,
\end{center}
Now, we can write  (\ref{090106})  as follows 
\begin{center}
$\omega = (\gamma_{1}+a_{1}^{1}df_{1}^{2}+a_{2}^{1}df_{2}^{2}+\cdots+a_{k}^{1}df_{k}^{2})\wedge dx_{1} + \cdots + (\gamma_{k}+a_{1}^{k}df_{1}^{2}+a_{2}^{k}df_{2}^{2}+\cdots++a_{k}^{k}df_{k}^{2})\wedge dx_{k} + (\mu_{1}+df_{1}^{2})\wedge dy_{1} + \cdots (\mu_{k}+df_{k}^{2})\wedge dy_{k}.$
\end{center}
Thus, we obtain the following expression for $\omega$:
$$\omega = \zeta +(-a_{1}^{1}dx_{1} - \cdots -a_{1}^{k}dx_{k} - dy_{1})\wedge df_{1}^{2}+\cdots+(-a_{k}^{1}dx_{1} - \cdots -a_{k}^{k}dx_{k} - dy_{k})\wedge df_{k}^{2},$$
where $\zeta\in \langle dx_{1},\dots,dx_{k},dy_{1},\dots,dy_{k}\rangle$ is a linear 2-form. This is,
\begin{equation*}
\omega = \zeta + dt_{1}\wedge df_{1}^{2} + \cdots+ dt_{k}\wedge df_{k}^{2},
\end{equation*}
with $t_{1},\dots,  t_{k}$ linear and dependent only on the variables $(x_{1},\dots, x_{k}, y_{1}, \dots, y_{k})$. Note that
\begin{equation*}
d\omega = d\zeta,
\end{equation*}
and since  $0=d\omega = d\zeta$ and $\zeta\in \langle dx_{1},\dots,dx_{k},dy_{1},\dots,dy_{k}\rangle$  we conclude that $\zeta$ is a linear $2$-form  which  dependents only on the variables $(x_{1},\dots, x_{k}, y_{1}, \dots, y_{k})$.

\begin{flushright}
    $\Box$
\end{flushright}

\section{Proof of Theorem \ref{classif}}
In order to proof this Theorem  we will use a similar idea in  \cite{LPT}. Indeed,  we will use  the  Theorem \ref{TeoDarb} and Jouanolou's lemma.

Let $\theta$  be  the homogeneous polynomial $1$-form in $\mathbb{C}^{n+1}$ of degree $2$  and of class $k$  which induces $\mathcal{F}$. Consider the polynomial $2$-forma $ d\theta.$  Since   the classe  of  $\theta$   is $k$ we have  that    $(d\theta)^{k+1}\not\equiv0$ and  $(d\theta)^{k+2} \equiv0$. It follows from  Theorem \ref{TeoDarb} that  there is a  coordinates system  $(x_{0},\dots,x_{k},y_{0},\dots,y_{k},z)$  with $z=\{z_{2k+3},\dots,z_{n+1}\}$ such $d\theta$ 
reduces to one of the following normal forms:
\begin{enumerate}
\item 
$
d\theta = \pi^*\eta, 
$
where  $\eta$  is a   closed  homogeneous $2$-form  of degree one  on  $\mathbb{C}^{2k+2}$  and   $\pi(x_{0},\dots, x_{k}, y_{0}, \dots, y_{k},z)=(x_{0},\dots, x_{k}, y_{0}, \dots, y_{k}) $ denotes the canonical linear projection. 
\item
$d\theta=   \pi^* \vartheta +dt_{1}\wedge dh_{1} +\cdots +dt_{k}\wedge dh_{k},$
where  $\vartheta $  is a   closed  homogeneous $2$-form  of degree one  on  $\mathbb{C}^{2k+2}$, $t_{1},\dots, t_{k}$ are linear  polynomials  in the variables $(x_{0},\dots, x_{k}, y_{0}, \dots, y_{k})$ and $h_{1}, \dots, h_{k}$ are quadratic  polynomials.
\end{enumerate}

Since $\theta$  induces a distribution on $\mathbb{P}^n$, then 
$i_{R}\theta = 0$ and  it follows from Jouanolou's lemma (\ref{170101}) that 
\begin{equation}\label{180103}
i_{R}d\theta=3\theta.
\end{equation}
If we write $R= R_1+R_2$, where $R_1=\sum_jx_j\frac{\partial}{\partial x_j} + y_j\frac{\partial}{\partial y_j}$ and  $R_2=\sum_{i}z_i \frac{\partial}{\partial z_i} $, it is clear that 
$i_{R}( \pi^* \eta)=\pi^* (i_{R_1}  \eta)$.
Then, 
$$3\theta= i_{R}( d \theta)= i_{R} (\pi^*\eta)=\pi^* (i_{R_1}  \eta).$$
This proves the case $i)$. 

Now, we have that 
\begin{equation}\label{180104}
i_{R}\sum_{i}(dt_{i}\wedge dh_{i})= \sum_{i}[(i_{R}dt_{i})dh_{i} - (i_{R}dh_{i})dt_{i}]=\sum_{i}(t_{i}dh_{i} - 2h_{i}dt_{i}).
\end{equation}
On the other hand, by contracting  $d\theta$   with the radial vector field  $R$ we obtain
\begin{equation}\label{180105}
3\theta= i_{R}d\theta = i_{R}( \pi^* \vartheta)+i_{R}(\sum_{i=1}^{k}dt_{i}\wedge dh_{i}).
\end{equation}

Substituting (\ref{180103}) and (\ref{180104}) in (\ref{180105}) we conclude that 
$$\theta = \frac{1}{3}\left[  \pi^* (i_{R_1}  \vartheta)+  \sum_{i}(t_{i}dh_{i} - 2h_{i}dt_{i})\right].$$

Finally,  we  observe that  the form $  \sum_{i}(t_{i}dh_{i} - 2h_{i}dt_{i})$ is the pull-back of the canonical contact form
$$
\theta_0=  \sum_{i}(u_{i}dw_{i} - 2u_{i}dw_{i}) 
$$
via the rational map
$${\scriptsize \xi: (x_{0}:\cdots: x_{k}: y_{0}: \cdots: y_{k}:z)\in \mathbb{P}^n \dashrightarrow  (t_0:\dots:t_{k}: h_0,\dots:h_{k}) \in \mathbb{P}(1^{k+1},2^{k+1}).}$$

\begin{flushright}
    $\Box$
\end{flushright}

\end{document}